# COMPUTABLE EXPONENTIAL BOUNDS FOR SCREENED ESTIMATION AND SIMULATION

By Ioannis Kontoyiannis [1] and Sean P. Meyn [2]

*Athens University of Economics and Business and University of Illinois, Urbana-Champaign*

Suppose the expectation $E(F(X))$ is to be estimated by the empirical averages of the values of $F$ on independent and identically distributed samples $\{X_i\}$. A sampling rule called the "screened" estimator is introduced, and its performance is studied. When the mean $E(U(X))$ of a different function $U$ is known, the estimates are "screened," in that we only consider those which correspond to times when the empirical average of the $\{U(X_i)\}$ is sufficiently close to its known mean. As long as $U$ dominates $F$ appropriately, the screened estimates admit exponential error bounds, even when $F(X)$ is heavy-tailed. The main results are several nonasymptotic, explicit exponential bounds for the screened estimates. A geometric interpretation, in the spirit of Sanov's theorem, is given for the fact that the screened estimates always admit exponential error bounds, even if the standard estimates do not. And when they do, the screened estimates' error probability has a significantly better exponent. This implies that screening can be interpreted as a variance reduction technique. Our main mathematical tools come from large deviations techniques. The results are illustrated by a detailed simulation example.

**1. Introduction.** Suppose we wish to estimate the expectation,

$$\mu := E[X^{3/4}] = \int_1^\infty x^{3/4} f(x)\, dx,$$

based on $n$ independent samples $X_1, X_2, \ldots, X_n$ drawn from some unknown density $f$ on $[1, \infty)$. Suppose, also, we have reasons to suspect that $f$ has a fairly heavy right tail, and assume that the only specific piece of information

Received November 2006; revised July 2007.
[1]Supported in part by a Sloan Foundation Research Fellowship.
[2]Supported in part by NSF Grant ECS 99-72957.
*AMS 2000 subject classifications.* Primary 60C05, 60F10; secondary 60G05, 60E15.
*Key words and phrases.* Estimation, Monte Carlo, simulation, large deviations, computable bounds, measure concentration, variance reduction.

we have available is the value of the mean of $f$, $\nu := E(X) = \int_1^\infty x f(x)\, dx$, perhaps also its variance. Because of the heavy right tail, it is natural to expect significant variability in the data $\{X_i\}$ as well as in the subsequent estimates of $\mu$. For definiteness, assume that the unknown density is $f(x) = \frac{5}{2x^{7/2}}$, for $x \geq 1$ [and $f(x) = 0$, otherwise], so that $\mu = 10/7$ and $\nu = 5/3$.

Consider the simplest (and most commonly used) estimator for $\mu$; for each $k \leq n$, let $\hat{S}_k$ denote the empirical average of the transformed samples $\{X_i^{3/4}\}$,

$$\hat{S}_k := \frac{1}{k} \sum_{i=1}^{k} X_i^{3/4}, \qquad 1 \leq k \leq n.$$

Although the law of large numbers guarantees that the sequence of estimates $\{\hat{S}_k\}$ is consistent and the central limit theorem implies that the rate of convergence is of order $n^{-1/2}$, a quick glance at the behavior of $\hat{S}_k$ for finite $k$ reinforces the concern that the estimates are highly variable: The plots in Figure 1 clearly indicate that, up to $k = n = 5000$, the $\{\hat{S}_k\}$ are still quite far from having converged.

Since $f$ is heavy-tailed, this irregular behavior is hardly surprising: Indeed, as $n$ grows, the error probability $\Pr\{\hat{S}_n > \mu + \varepsilon\}$ decays like

(1.1) $$\Pr\{\hat{S}_n > \mu + \varepsilon\} \sim \frac{1}{\varepsilon^{10/3} n^{7/3}}, \qquad n \to \infty,$$

for any $\varepsilon > 0$; see, for example [12]. Therefore, unlike with most classical exponential error bounds, here the error probability decays polynomially in the sample size $n$, and with a rather small power at that.

This state of affairs is discouraging, but suppose we decide to use the additional information we have about $f$, namely that its mean $\nu$ equals $5/3$, in order to "screen" the estimates $\{\hat{S}_k\}$. This can be done as follows: Together with the $\{\hat{S}_k\}$, also compute the empirical averages $\{\hat{T}_k\}$ of the samples $\{X_i\}$ themselves,

$$\hat{T}_k = \frac{1}{k}\sum_{i=1}^{k} X_i, \qquad 1 \leq k \leq n,$$

and *only consider estimates* $\hat{S}_k$ *at times $k$ when the corresponding average $\hat{T}_k$ is within a fixed threshold $u > 0$ from its known mean.* That is, only examine $\hat{S}_k$ if at that same time $k$, $|\hat{T}_k - \nu| < u$.

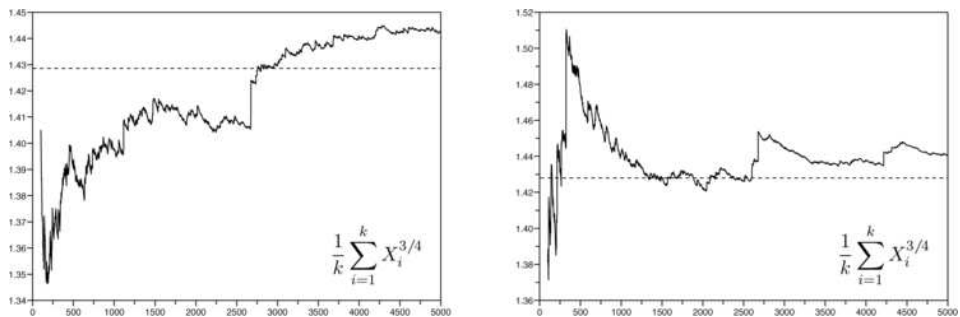

FIG. 1.  *Two typical realizations of the estimates $\{\hat{S}_k\}$ for $k = 100, 101, \ldots, n = 5000$.*

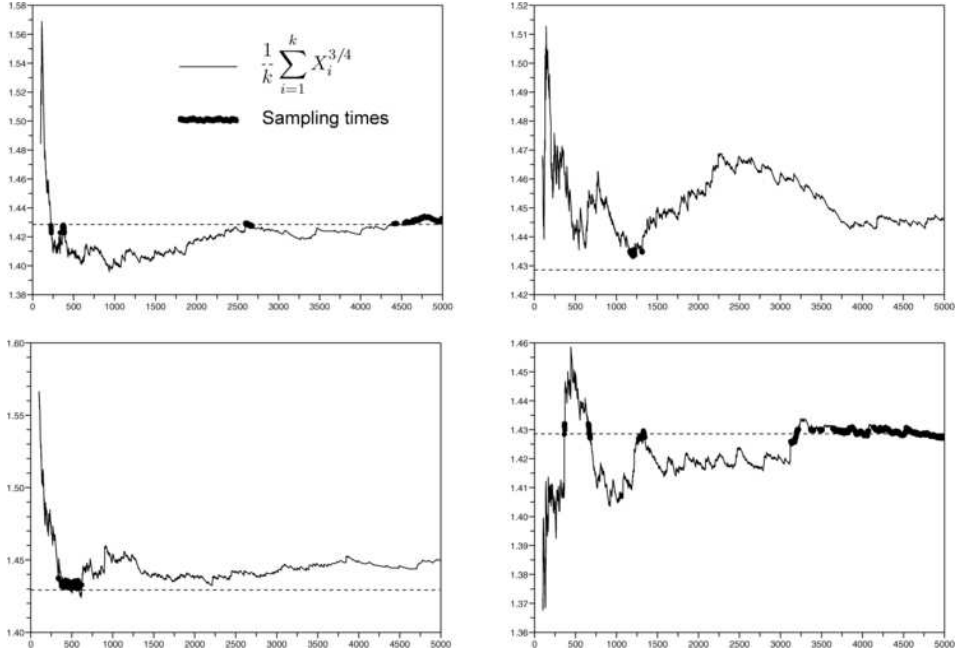

FIG. 2. *Four typical realizations of the estimates $\{\hat{S}_k\}$ for $k = 100, 101, \ldots, n = 5000$. The "screened estimates" are plotted in bold, and they are simply the original $\hat{S}_k$ at times $k$ when the corresponding empirical average $\hat{T}_k$ is within $u = 0.005$ of its mean $\nu = 5/3$.*

This results in what we call in this paper the "*screened estimator*" of $\mu$. Figure 2 illustrates its performance on four different realizations of the above experiment.

More generally, assume $X, X_1, X_2, \ldots$ are independent and identically distributed (i.i.d.) random variables with unknown distribution, and we wish to estimate the expectation $\mu := E[F(X)]$ for a given function $F : \mathbb{R} \to \mathbb{R}$, while we happen to know the value of the expectation $\nu := E[U(X)]$ of a different function $U : \mathbb{R} \to \mathbb{R}$. In this general setting, we introduce:

**The Screened Estimator**. For each $k \geq 1$, together with the empirical averages $\{\hat{S}_k\}$ of the $\{F(X_i)\}$ also compute the averages $\{\hat{T}_k\}$ of the $\{U(X_i)\}$, and *only* consider estimates $\hat{S}_k$ at times $k$ when $\hat{T}_k$ is within a fixed threshold $u > 0$ from its mean, that is, $|\hat{T}_k - \nu| < u$.

The intuition is simple. In cases when we suspect that the empirical distribution $\hat{P}_k$ of the samples $\{X_i : i \leq k\}$ is likely to be far from the true underlying distribution $P$, we can check that the projection $\int U \, d\hat{P}_k = \hat{T}_k$ of $\hat{P}_k$ along a function $U$ is close to the projection $\int U \, dP = \nu$ of the true distribution $P$ along $U$. Of course this does not guarantee that $\hat{P}_k \approx P$ or that $\hat{S}_k \approx \mu$, but it *does* rule out instances $k$ when it is certain that $\hat{P}_k$ differs significantly from $P$.

More importantly, as we shall see next, it is often possible to obtain *explicitly computable exponential error bounds* for the screened estimator, even when the error probability of the standard estimates $\{\hat{S}_k\}$ decays at a polynomial rate.

The purpose of this paper is twofold. First, we provide a theoretical explanation for the practical advantage of the screened estimator: We develop general conditions under which the error probability of the screened estimator decays exponentially, regardless of the tail of the distribution of the $\{F(X_i)\}$. The main assumption is that $U$ dominates $F$ from above, in that $\sup_x [F(x) - \beta U(x)]$ is finite for all $\beta > 0$, where the supremum is over all $x$ in the support of $X$. Then we state and prove a number of explicit exponential bounds for the error probability of the screened estimator, which are easily computable and readily applicable to specific problems where the only information we have about the unknown underlying distribution is the mean and perhaps also the variance of $U(X)$ for a particular function $U$.

To illustrate, we return to the example of estimating the expectation $\mu = E(X^{3/4})$ with respect to an unknown density $f$ on $[1, \infty)$, based on $n$ i.i.d. samples $X_1, \ldots, X_n$ drawn from $f$, and assuming that we only know the mean (and perhaps some higher moments) of $X$. In the above notation, this corresponds to $F(x) \equiv x^{3/4}$ and $U(x) \equiv x$. The proof of the following proposition is given at the end of Section 3.

PROPOSITION 1.1. (i) *The error probability of the standard estimator* $\{\hat{S}_n\}$ *decays to zero at a polynomial rate: If the density $f$ is given by $f(x) = \frac{5}{2x^{7/2}}$ for $x \geq 1$, then for any $\varepsilon > 0$,*

$$\Pr\{\hat{S}_n - \mu > \varepsilon\} \sim \frac{1}{\varepsilon^{10/3} n^{7/3}}, \qquad n \to \infty.$$

(ii) *The error probability of the screened estimator decays to zero exponentially fast: If the only information we have about $f$ is that its mean $\nu$ equals $5/3$, then we can conclude that for all $\varepsilon, u > 0$ there exists $I(\varepsilon, u) > 0$ such that*

$$\Pr\{\hat{S}_n - \mu > \varepsilon \text{ and } |\hat{T}_n - \tfrac{5}{3}| < u\} \leq e^{-nI(\varepsilon, u)} \qquad \text{for all } n \geq 1.$$

SCREENED ESTIMATION AND SIMULATION 5(iii) *If, in addition, we know that the variance of $f$ equals $20/9$, then an explicit exponential bound can be computed: For any $\varepsilon > 0$ and any $0 < u \leq \frac{\varepsilon}{20}$,*

$$\Pr\{\hat{S}_n - \mu > \varepsilon \text{ and } |\hat{T}_n - \tfrac{5}{3}| < u\} \leq e^{-(0.005) \times n\varepsilon^2} \qquad \text{for all } n \geq 1.$$

(iv) *If we also know that the value of the covariance between $X^{3/4}$ and $X$ under $f$ is $20/21$, then the following more accurate bound can be obtained: For any $\varepsilon > 0$ and any $0 < u \leq \frac{\varepsilon}{20}$,*

(1.2) $$\Pr\{\hat{S}_n - \mu > \varepsilon \text{ and } |\hat{T}_n - \tfrac{5}{3}| < u\} \leq e^{-(0.0367) \times n\varepsilon^2},$$

*for all $n \geq 1$.*

As long as the mean of $X$ is known, we can employ the screened estimator and be certain that it will have an exponentially small error probability, whereas the standard estimator's probability of error may decay at least as slowly as $n^{-7/3}$. If the variance of $X$ is also known, then for the specific values in the simulation examples in Figure 2, with $\varepsilon = 0.2$, $u = 0.005$ and $n = 5000$, part (iii) of the proposition gives

$$\Pr\{\hat{S}_n - \mu > 0.2 \text{ and } |\hat{T}_n - \tfrac{5}{3}| < 0.005\} \leq 0.368.$$

This is fairly weak, despite the fact that $\varepsilon = 0.2$ is a rather moderate margin of error. But the error probability does decay exponentially, and with $n = 10{,}000$ samples the corresponding upper bound is only $\approx 0.136$, while for $n = 15{,}000$ it is $\approx 0.0498$. And if, in addition, the value of the covariance between $X^{3/4}$ and $X$ is available, then part (iv) gives a much more accurate result even for smaller $\varepsilon$: Taking $\varepsilon = 0.1$, $u = 0.005$ and $n = 5000$,

$$\Pr\{\hat{S}_n - \mu > 0.1 \text{ and } |\hat{T}_n - \tfrac{5}{3}| < 0.005\} \leq 0.1596,$$

and for $n = 10{,}000$ samples the corresponding bound is $\approx 0.025$.

Two points of caution are in order here. The first is perhaps somewhat subtle and has to do with the interpretation of the above error bounds. What exactly does (1.2) say? Is it the case that, at any time $k$ when $\hat{T}_k$ is within $u$ of its mean, we can apply (1.2) to obtain a bound on the probability of error for the corresponding estimate $\hat{S}_k$? Strictly speaking, the answer is "no"; since the times at which the screening averages $\{\hat{T}_k\}$ are close to their mean are random, (1.2) cannot be automatically invoked. A strict operational interpretation of the mathematical statement in (1.2) is as follows: First choose and fix an $n$ such that (1.2) offers a satisfactory guarantee on the error probability; here $n$ may be the total number of samples available, or it may be the number of samples we decide to generate from $f$. Then look at $\hat{T}_n$, and if $|\hat{T}_n - \nu| < u$, it is legitimate to use the error bound (1.2) for the value of the estimate $\hat{S}_n$ at the last sample time $n$. Otherwise, do not use the bound (1.2) at all.

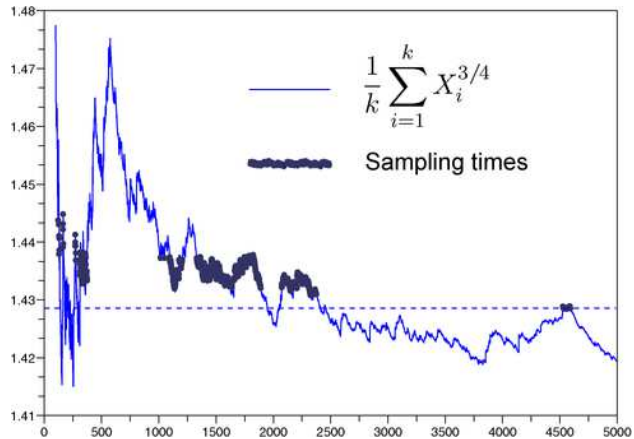

FIG. 3. *Another realization of the empirical estimates $\{\hat{S}_k\}$ for $k = 100, 101, \ldots, n = 5000$, plotted together with the screened estimates shown in bold (where $u = 0.005$ as before). The screened estimates at earlier times are less accurate than some of the later estimates that are ignored by the screened estimator.*

The same interpretation applies to any application of the screened estimator. On the one hand, screening gives a powerful heuristic for selecting times $k$ when the $\hat{S}_k$ are more likely to be accurate, and it can be used as a diagnostic tool to actually *rule out* times $k$ when it is certain that the empirical distribution of the samples is not close to the true underlying distribution. On the other hand, in cases when it is required that the error probability be precisely quantified, the sampling times cannot be random and they have to be decided upon in advance.

The second point is based on some results we observed in simulation experiments, indicating that the sampling times $k$ picked out by the screened estimator are not all equally reliable: Naturally, since the probability of error decays exponentially, earlier times correspond to much looser error bounds, while the error probability of estimates obtained during later times can be more tightly controlled. This is illustrated by the (rather atypical but not impossibly rare) results shown in Figure 3.

From the probabilistic point of view, the following calculation gives a quick explanation for the fact that the screened estimator leads to exponential error bounds in great generality (although this is not how the actual error bounds in Section 3 are obtained). Suppose the $\{\hat{S}_k\}$ are used to estimate the mean $\mu = E(F(X))$ for some $F$, while we know $\nu = E(U(X))$ for a different function $U$ that dominates $F$ in that $\operatorname{ess\,sup}_X[F(X) - \beta U(X)] < \infty$, for all $\beta > 0$. Although $F(X)$ may be heavy-tailed, in which case the $\{\hat{S}_k\}$ themselves will not admit exponential error bounds, the error probability of the screened estimator is bounded by

$$\Pr\{\hat{S}_n - \mu > \varepsilon \text{ and } |\hat{T}_n - \nu| < u\}$$
(1.3)
$$\leq \Pr\left\{\frac{1}{n}\sum_{i=1}^{n}[F(X_i) - \beta U(X_i)] - (\mu - \beta\nu) > \varepsilon - \beta u\right\}.$$

Since $E[F(X) - \beta U(X)] = \mu - \beta\nu$, for $0 < \beta < \frac{\varepsilon}{u}$ this is a large deviations probability for the right tail of the partial sums of the random variables $\{F(X_i) - \beta U(X_i)\}$, which are (a.s.) *bounded above*. It is, therefore, no surprise that this probability is exponentially small.

1.1. *Screening and control variates.* A well-known and commonly used technique for reducing the variance of an estimator in classical Monte Carlo simulation is the method of *control variates*; see, for example, the standard texts [7, 11, 13] or the paper [9] for extensive discussions. This method is based on the observation that in many applications—exactly as in our setting—there is a function $U$ whose expectation $\nu = E[U(X)]$ is known. Therefore, replacing the estimates $\{\hat{S}_k\}$ for $\mu = E[F(X)]$ with the *control variate estimates*,

$$\tilde{S}_k := \frac{1}{k}\sum_{i=1}^{k}(F(X_i) - \beta[U(X_i) - \nu]), \qquad 1 \leq k \leq n,$$

yields an estimator which is still consistent (since the additional term has zero mean) but whose variance is different from that of $\{\hat{S}_k\}$. In fact, choosing (or estimating) the value of the constant $\beta$ appropriately always leads to an estimator with strictly reduced variance, as long as $F(X)$ and $U(X)$ are correlated random variables.

This technique is widely employed in practice; see the references above as well as [1, 6]; also the text [8] contains many examples of current interest in computational finance and pointers to the relevant literature. In particular, functions $U$ that appear in applications as control variates provide a natural class of screening functions that can be incorporated in the design on the screened estimator.

An interesting connection between these two methods (control variates and screening) is seen in that the second probability in (1.3) above is exactly the error probability for the control variate estimates $\{\tilde{S}_k\}$. More generally, in cases where control variates (or some other method) are used to reduce the variance of the $\{\hat{S}_k\}$, we view screening as a sampling rule which complements (and does not replace) variance reduction or variance estimation techniques. The connection between screening and variance reduction is an intriguing one, and will be explored in subsequent work [10].

1.2. *Outline and summary of results.* The general results in Sections 2 and 3 parallel those presented for the example in Proposition 1.1. Theorems 2.1, 2.2 and 2.3 offer a theoretical description of the large deviations behavior of the screened estimator's error probability, both asymptotically and for finite $n$. The only assumptions necessary are that $E(F(X))$ is finite, and that the mean $E(U(X))$ is known for some function $U$ which dominates $F$ in that $\operatorname{ess\,sup}[F(X) - \beta U(X)] < \infty$ for all $\beta > 0$. Then the error probability admits a nontrivial exponential bound, regardless of the distribution of $X$. The exponent can be expressed either as a Fenchel–Legendre transform or in terms of relative entropy, and the relative entropy formulation leads to an elegant geometric explanation for the fact that the screened estimator's error probability always decays exponentially.

When $F(X)$ and $U(X)$ also have finite second moments, and assuming that the variance $\operatorname{Var}(U(X))$ is known, in Theorem 3.1 we give an explicit, easily computable, exponential bound for the error probability. The bound holds for all $n \geq 1$, and the exponent is of order $\varepsilon^2$ for small $\varepsilon, u$. Also, a more refined bound is given when the value of the covariance between $F(X)$ and $U(X)$ is available. These are the main results of this paper.

In Section 4 we consider the case when $F(X)$ and $U(X)$ have finite exponential moments, so that the standard estimator $\{\hat{S}_k\}$ already has an exponentially vanishing error probability. Theorem 4.1 shows that the screened estimator's error probability decays at a strictly faster exponential rate, and the difference between the exponents is more precisely quantified in Theorem 4.2: It is shown to be of order $\varepsilon^2$ for small $\varepsilon, u$, and this is used to draw a different heuristic connection between screening and variance reduction techniques.

Section 5 contains the proofs of Theorems 2.1, 2.2 and 2.3.

Finally, we mention that the screening idea can also be applied in the context of more complex problems arising in Markov chain Monte Carlo (MCMC) simulation. Such generalizations are by no means immediate, and they will be explored in subsequent work.

**2. Large deviations.** In this section we give a theoretical explanation for the (sometimes dramatic) performance improvement offered by the screened estimator. For explicit bounds like those presented in the Introduction, see Section 3.

Let $X, X_1, X_2, \ldots$ be i.i.d. random variables with common law given by the probability measure $P$ on $\mathbb{R}$. Given a function $F : \mathbb{R} \to \mathbb{R}$ whose mean is to be estimated by the empirical averages of the $\{F(X_i)\}$, for the purposes of this section only we consider a slightly simplified version of the screened estimator: Assuming the mean $\nu = E(U(X))$ of a different function $U : \mathbb{R} \to \mathbb{R}$ is known, we examine the screened estimator based on the one-sided screening event, $\{\sum_{i=1}^{n} U(X_i) - n\nu < nu\}$, for some $u > 0$. To avoid cumbersome notation, write $S_n := \sum_{i=1}^{n} F(X_i)$ and $T_n := \sum_{i=1}^{n} U(X_i)$, $n \geq 1$.

In the first result, Theorem 2.1 below, we obtain representations for the asymptotic exponents of the error probability, both for the standard estimator and for the screened estimator. The exponents are expressed in terms of relative entropy, in the spirit of Sanov's theorem; see [3, 4, 14]. Recall that the relative entropy between two probability measures $P$ and $Q$ on the same space is defined by

$$H(P\|Q) := \begin{cases} \int dP \log \frac{dP}{dQ}, & \text{when } \frac{dP}{dQ} \text{ exists,} \\ \infty, & \text{otherwise.} \end{cases}$$

Theorem 2.1 follows from the more general results in Theorems 2.2 and 2.3 below; its proof is given in Section 5.

THEOREM 2.1 (Sanov asymptotics). *Suppose the functions $F : \mathbb{R} \to [0, \infty)$ and $U : \mathbb{R} \to \mathbb{R}$ have finite first moments $\mu := E[F(X)]$, $\nu := E[U(X)]$, and also finite second moments, $E[F(X)^2]$, $E[U(X)^2]$. Assume that $F(X)$ is heavy-tailed in that $E[e^{\theta F(X)}] = \infty$ for all $\theta > 0$, and that $U$ dominates $F$ in that $m(\beta) := \operatorname{ess\,sup}[F(X) - \beta U(X)] < \infty$ for all $\beta > 0$. Then:*

(i) *The error probability of the standard estimator decays subexponentially: For all $\varepsilon > 0$,*

$$\lim_{n \to \infty} \frac{1}{n} \log \Pr\{S_n - n\mu > n\varepsilon\} = -\inf_{Q \in \Sigma} H(Q\|P) = 0,$$

*where $\Sigma$ is the set of all probability measures $Q$ on $\mathbb{R}$ such that $\int F dQ - \mu > \varepsilon$.*

(ii) *The error probability of the screened estimator decays exponentially: For all $\varepsilon, u > 0$,*

$$\lim_{n \to \infty} \frac{1}{n} \log \Pr\{S_n - n\mu > n\varepsilon \text{ and } T_n - n\nu < nu\} = -\inf_{Q \in E} H(Q\|P) < 0,$$

*where $E \subset \Sigma$ is the set of all probability measures $Q$ on $\mathbb{R}$ such that $\int F \, dQ - \mu > \varepsilon$ and $\int U \, dQ - \nu < u$.*

Therefore, while the (asymptotic) exponent of the error probability of the standard estimator is equal to zero, the exponent of the error probability of the screened estimator is strictly positive. Although this situation is only possible when the relative entropy is minimized over an infinite-dimensional space of measures [in that the exponent $\inf_{Q \in \Sigma} H(Q\|P)$ cannot be zero when $X$ takes only finitely many values], it is perhaps illuminating to offer a geometric description.

The large oval in the first diagram in Figure 4 depicts the space of all probability measures $Q$ on $\mathbb{R}$, and the small "cap" on the left is the set $\Sigma$ of those $Q$ with $\int F dQ - \mu > \varepsilon$. The gray shaded area corresponds to the "smallest" subset of $\Sigma$ such that the infimum of $H(Q\|P)$ over this subset is zero. (Of course this set is not exactly well defined, but it does convey the correct intuition.) In the second diagram, the black shaded area corresponds to set $E$, formed by the intersection of $\Sigma$ with the half-space $H = \{Q : \int U\,dQ - \nu < u\}$. Note that $H$ is a "typical" set under $P$, in that $P \in H$ and the empirical measure of the $\{X_i\}$ will eventually concentrate there by the law of large numbers. Nevertheless, when $\Sigma$ is intersected with $H$ to give $E$, Theorem 2.1 tells us that it excludes the part of $\Sigma$ which is close to $P$ in relative entropy (the gray area), and this forces the result of the minimization over $Q \in E$ to be strictly positive; the limiting minimizer $Q^*$, assuming it exists, is shown as laying on the common boundary of $\Sigma$ and $H$.

The following two theorems give a more precise and complete description of the large deviations properties of the probabilities of interest. Formally, they simply establish a version of Cramér's theorem in the present setting. What is perhaps somewhat surprising is that this is done without *any* assumption of finite exponential moments. In the presence of the domination condition $m(\beta) < \infty$, it turns out that it is only necessary to assume finite first (and in some cases second) moments for $F(X)$ and $U(X)$.

The results in Theorems 2.2 and 2.3 will form the basis for the development of the bounds in Section 3. Their proofs are given in Section 5.

THEOREM 2.2 (Exponential upper bounds). *Suppose the functions* $F : \mathbb{R} \to \mathbb{R}$ *and* $U : \mathbb{R} \to \mathbb{R}$ *are such that* $\mu := E[F(X)]$ *and* $\nu := E[U(X)]$ *are both finite, and that* $m(\beta) := \operatorname{ess\,sup}[F(X) - \beta U(X)] < \infty$ *for all* $\beta > 0$. *Then for all* $\varepsilon, u > 0$:

(i) $\Pr\{S_n - n\mu > n\varepsilon, T_n - n\nu < nu\} \leq \exp\{-nH(E\|P)\}$, *for all* $n \geq 1$, *where*,

(2.1) $$H(E\|P) := \inf\{H(Q\|P) : Q \in E\},$$

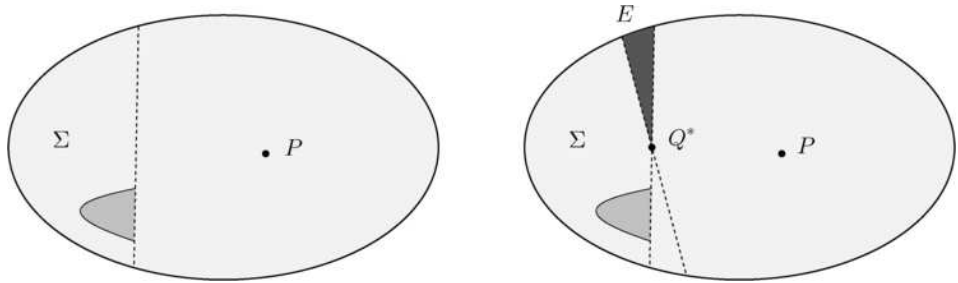

FIG. 4. *Geometric illustration of the fact that* $\inf_{Q \in \Sigma} H(Q\|P) = 0$ *whereas* $\inf_{Q \in E} H(Q\|P)$ *is strictly positive.*

and $E$ is the set of all probability measures $Q$ on $\mathbb{R}$ such that $\int F\,dQ - \mu > \varepsilon$ and $\int U\,dQ - \nu < u$.

(ii) $\Pr\{S_n - n\mu > n\varepsilon, T_n - n\nu < nu\} \leq \exp\{-n\Lambda_+^*(\varepsilon, u)\}$, for all $n \geq 1$, where,

$$\Lambda_+^*(\varepsilon, u) := \sup_{\theta_1, \theta_2 \geq 0} \{\theta_1(\mu + \varepsilon) - \theta_2(\nu + u) - \Lambda_+(\theta_1, \theta_2)\},$$

with $\Lambda_+(\theta_1, \theta_2) := \log E[\exp\{\theta_1 F(X) - \theta_2 U(X)\}]$, $\theta_1, \theta_2 \geq 0$.

(iii) The rate function $\Lambda_+^*(\varepsilon, u)$ is strictly positive.

THEOREM 2.3 (Large Deviations Asymptotics). *Under the assumptions of Theorem 2.2, if, in addition, $F(X)$ and $U(X)$ have finite second moments, then for all $\varepsilon, u > 0$,*

$$(2.2) \qquad \lim_{n \to \infty} \frac{1}{n} \log \Pr\{S_n - n\mu > n\varepsilon, T_n - n\nu < nu\} = -\Lambda_+^*(\varepsilon, u),$$

*and $\Lambda_+^*(\varepsilon, u)$ coincides with the rate function $H(E\|P)$ given in (2.1).*

**3. Bounds for arbitrary tails.** Let $X, X_1, X_2, \ldots$ be i.i.d. random variables. Given functions $F, U : \mathbb{R} \to \mathbb{R}$, write $S_n = \sum_{i=1}^n F(X_i)$ and $T_n = \sum_{i=1}^n U(X_i)$. We begin by restating part of Theorem 2.2. Since the two-sided error event $\{S_n - n\mu > n\varepsilon, |T_n - n\nu| < nu\}$ is contained in $\{S_n - n\mu > n\varepsilon, T_n - n\nu < nu\}$, we have:

COROLLARY 3.1. *Suppose the functions $F : \mathbb{R} \to \mathbb{R}$ and $U : \mathbb{R} \to \mathbb{R}$ are such that $\mu := E[F(X)]$ and $\nu := E[U(X)]$ are both finite, and that $m(\beta) := \operatorname{ess\,sup}[F(X) - \beta U(X)] < \infty$ for all $\beta > 0$. Then for all $n \geq 1$ and all $\varepsilon, u > 0$,*

$$\Pr\{S_n - n\mu > n\varepsilon, |T_n - n\nu| < nu\} \leq e^{-n\Lambda_+^*(\varepsilon, u)},$$

*where the exponent, $\Lambda_+^*(\varepsilon, u)$, is given by*

$$\sup_{\theta_1 \geq 0, \theta_2 \geq 0} \{\theta_1(\mu + \varepsilon) - \theta_2(\nu + u) - \log E[\exp\{\theta_1 F(X) - \theta_2 U(X)\}]\},$$

*and is strictly positive.*

REMARKS. 1. An exactly analogous result holds if instead of $m(\beta)$ we assume that $\operatorname{ess\,sup}[F(X) + \beta U(X)] < \infty$ for all $\beta > 0$. Then, repeating the Chernoff argument in the proof of Theorem 2.2 for the one-sided error event $\{S_n - n\mu > \varepsilon, T_n - n\nu > -nu\}$ leads to the same bound, but with the exponent, $\Gamma_+^*(\varepsilon, u)$, given by

$$\sup_{\theta_1 \geq 0, \theta_2 \geq 0} \{\theta_1(\mu + \varepsilon) + \theta_2(\nu - u) - \log E[\exp\{\theta_1 F(X) + \theta_2 U(X)\}]\},$$

and $\Gamma_+^*(\varepsilon, u)$ can be similarly seen to be strictly positive.

2. Replacing $F$ by $-F$ yields a corresponding result for the left tail. If $\operatorname{ess\,inf}[F(X) + \beta U(X)] > -\infty$ for all $\beta > 0$,

$$\Pr\{S_n - n\mu < -n\varepsilon, |T_n - n\nu| < nu\} \leq e^{-n\Lambda_-^*(\varepsilon,u)},$$

where $\Lambda_-^*(\varepsilon, u)$ is given by

$$\sup_{\theta_1 \geq 0, \theta_2 \geq 0} \{\theta_1(-\mu + \varepsilon) + \theta_2(-\nu - u) - \log E[\exp\{-\theta_1 F(X) - \theta_2 U(X)\}]\},$$

and is strictly positive. Moreover, in view of the previous remark, an analogous bound holds under the assumption that $\operatorname{ess\,inf}[F(X) - \beta U(X)] > -\infty$ for all $\beta > 0$; in this case the exponent is replaced by

$$\Gamma_-^*(\varepsilon, u) = \sup_{\theta_1 \geq 0, \theta_2 \geq 0} \{\theta_1(-\mu + \varepsilon) + \theta_2(\nu - u) - \log E[\exp\{-\theta_1 F(X) + \theta_2 U(X)\}]\},$$

which is also strictly positive.

3. Combining the observations in Remarks 1 and 2 immediately yields a bound on the two-sided deviations of $\{S_n\}$. If both $\mu = E[F(X)]$ and $\nu = E[U(X)]$ are finite, and also both $\operatorname{ess\,sup}[F(X) - \beta U(X)] < \infty$ and $\operatorname{ess\,inf}[F(X) + \beta U(X)] > -\infty$, for all $\beta > 0$, then for all $n \geq 1$ and all $\varepsilon, u > 0$,

$$\Pr\{|S_n - n\mu| > n\varepsilon, |T_n - n\nu| < nu\}$$
(3.1)
$$\leq e^{-n\Lambda_+^*} + e^{-n\Lambda_-^*} \leq 2e^{-n\min\{\Lambda_+^*, \Lambda_-^*\}},$$

where $\Lambda_+^*$ and $\Lambda_-^*$ are strictly positive. Although this double domination assumption may appear severe, it is generally quite easy to find functions $U$ that will dominate a given $F$ appropriately. For example, if $F(x) \equiv x$ we can simply take $U(x) \equiv x^2$, or, more generally, $U(x) \equiv x^{2k}$ for any positive integer $k$, assuming appropriately high moments exist.

4. In Remarks 1 and 2, two different domination assumptions were shown to give a bound on the right tail of the partial sums of $F$, and two more assumptions do the same for the left tail. Any of their four different combinations gives a bound similar to (3.1), with the appropriate combination of exponents.

If $F$ and $U$ also have finite second moments, an easily applicable, quantitative version of Corollary 3.1 can be obtained. The gist of the argument is the use of the boundedness of $[F(X) - \beta U(X)]$ in order to compute an explicit lower bound for the exponent $\Lambda_+^*(\varepsilon, u)$.

THEOREM 3.1. *Suppose* $E[F(X)] = E[U(X)] = 0$, $\mathrm{Var}(F(X)) \leq 1$, $\mathrm{Var}(U(X)) = 1$, *and that* $m(\beta) := \operatorname{ess\,sup}[F(X) - \beta U(X)] < \infty$ *for all* $\beta > 0$. *Then the following hold for all* $n \geq 1$:

(i) *For any* $\varepsilon, u > 0$, *if there exists* $\beta > 0$ *such that* $m(\beta) \leq \varepsilon - \beta u$, *then*
$$\Pr\{S_n > n\varepsilon,\ |T_n| < nu\} = 0.$$

(ii) *For any* $\varepsilon, u > 0$,

(3.2)
$$\log \Pr\{S_n > n\varepsilon, |T_n| < nu\}$$
$$\leq -2n \sup_{\alpha \in (0,1)} \left[\frac{m \cdot (1-\alpha)}{m^2 + 1 + (\alpha\varepsilon/u)^2 - 2\alpha\gamma\varepsilon/u}\right]^2 \varepsilon^2,$$

*where* $m := m(\frac{\alpha\varepsilon}{u})$ *and* $\gamma := E[F(X)U(X)]$ *is the covariance between* $F(X)$ *and* $U(X)$.

(iii) *Let* $K > 0$ *be arbitrary. Then for any* $\varepsilon > 0$ *and any* $0 < u \leq K\varepsilon$,

(3.3) $$\log \Pr\{S_n > n\varepsilon, |T_n| < nu\} \leq -\frac{n}{2}\left[\frac{M}{M^2 + (1 + 1/(2K))^2}\right]^2 \varepsilon^2,$$

*where* $M = m(\frac{1}{2K})$.

REMARKS. 1. The assumption that $\mathrm{Var}(F(X)) \leq 1$ in Theorem 3.1 seems to require that we know an upper bound on the variance of $F$ in advance, but in practice this is easily circumvented. In specific applications, we typically have a function $U$ that dominates $F$ in that, not only $m(\beta) < \infty$ for all $\beta > 0$, but also there are finite constants $C_1, C_2$ such that

(3.4) $\quad |F(x)| \leq C_1 U(x) + C_2 \quad$ for all $x$ in the support of $X$.

This is certainly the case for the example presented in the Introduction, as well as in the examples in Remark 3 above. A bound on the variance of $F(X)$ is obtained from (3.4), $\mathrm{Var}(F(X)) \leq C_1^2 \mathrm{Var}(U(X)) + C_2^2$. This and several other issues arising in the application on Theorem 3.1 are illustrated in detail in the proof of Proposition 1.1.

2. As will become clear from its proof, to use the bounds in Theorem 3.1 it is not necessary to know $m(\beta)$ exactly; any upper bound on the $\operatorname{ess\,sup}[F(X) - \beta U(X)]$ can be used in place of $m(\beta)$.

3. When $F(x) \equiv \tilde{F}(x) - \mu$, where $\mu$ is the unknown mean to be estimated, it is hard to imagine that the exact value of the covariance $\gamma$ may be known without knowing $\mu$. But, similarly to $m(\beta)$, in order to apply (3.2) it suffices to have an upper bound on $\gamma$, and such estimates are often easy to obtain. See the proof of Proposition 1.1 for an illustration.

4. The main difference between the bounds in (3.2) and (3.3) is that (3.3) only requires knowledge of the first and second moments of $U(X)$, whereas (3.2) also depends on $\gamma$. The bound in (3.3) is attractive because it is simple and it clearly shows that the exponent is of order $\varepsilon^2$ for small $\varepsilon$. Its main disadvantage is that it often leads to rather conservative estimates, since it ignores the potential correlation between $F(X)$ and $U(X)$ and it follows from (3.2) by an arbitrary choice for the parameter $\alpha$. The exponent in (3.2), on the other hand, despite its perhaps somewhat daunting appearance, is often easy to estimate and it typically gives significantly better results. This too is clearly illustrated by the results (and the proof) of Proposition 1.1.

5. Considering $-F$ in place of $F$ gives corresponding bounds for the lower tail of the partial sums $S_n$, under the assumption that $\operatorname{ess\,inf}[F(X)+\beta U(X)]$ be finite for all $\beta > 0$. As in (3.1), these can be combined with the corresponding results in Theorem 3.1 to give explicit exponential bounds for the two-sided deviation event, $\{|S_n| > \varepsilon, |T_n| < u\}$.

PROOF OF THEOREM 3.1. As already noted in (1.3) in the Introduction, for any $\beta > 0$,

$$\Pr\{S_n > n\varepsilon, |T_n| < u\} \le \Pr\left\{\frac{1}{n}\sum_{i=1}^{n}[F(X_i) - \beta U(X_i)] > \varepsilon - \beta u\right\}.$$

If the essential supremum $m(\beta)$ of the random variables $[F(X_i) - \beta U(X_i)]$ which are being averaged is smaller than the threshold $\varepsilon - \beta u$, then the above event is empty and its probability is zero, establishing (i).

Recall the definitions of $\Lambda_+$ and $\Lambda_+^*$ in Theorem 2.2. With any $\alpha \in (0,1)$, taking $\theta_2 = \alpha\varepsilon\theta_1/u$ in the definition of $\Lambda_+^*(\varepsilon, u)$, Corollary 3.1 yields

$$(3.5) \quad \Pr\{S_n > n\varepsilon, |T_n| < nu\} \le \exp\left\{-n\sup_{\theta \ge 0}[\theta(1-\alpha)\varepsilon - \Lambda_0(\theta)]\right\},$$

where $\Lambda_0(\theta) := \Lambda_+(\theta, \frac{\alpha\varepsilon\theta}{u})$. Write $s^2 := \operatorname{Var}(F(X)) \le 1$, define the random variable $Y := F(X) - \frac{\alpha\varepsilon}{u}U(X)$, and note that $Y \le m := m(\frac{\alpha\varepsilon}{u})$ a.s., $E(Y) = 0$, and

$$\operatorname{Var}(Y) = s^2 + \left(\frac{\alpha\varepsilon}{u}\right)^2 - \frac{2\alpha\varepsilon\gamma}{u} \le \sigma^2 := 1 + \left(\frac{\alpha\varepsilon}{u}\right)^2 - \frac{2\alpha\varepsilon\gamma}{u}.$$

Throughout the rest of the proof we assume, without loss of generality, that $m > 0$. [We know $m \ge 0$ by our assumptions, so if $m = 0$, then $\Lambda_0(\theta) \equiv 0$ and the supremum in (3.5) equals $+\infty$, implying that the probability of interest equals zero and that all the bounds stated in the theorem are trivially valid.]

Now we apply Bennett's lemma [4], Lemma 2.4.1, to get an upper bound on $\Lambda_0(\theta)$ as

$$\Lambda_0(\theta) = \log E(e^{\theta Y}) \le \log\left\{\frac{m^2}{m^2 + \sigma^2}e^{-\theta\sigma^2/m} + \frac{\sigma^2}{m^2 + \sigma^2}e^{\theta m}\right\}.$$

Using this and replacing $\theta$ by $\lambda/m$, the supremum in (3.5) is bounded below by

$$I(\varepsilon,\alpha,u) := \sup_{\lambda \geq 0}\left[\lambda x - \log\left\{\frac{e^{-\lambda\tau^2} + \tau^2 e^{\lambda}}{1+\tau^2}\right\}\right],$$

where

$$x := \frac{(1-\alpha)\varepsilon}{m} \quad \text{and} \quad \tau^2 := \frac{\sigma^2}{m^2} = \frac{1+((\alpha\varepsilon)/u)^2 - (2\alpha\varepsilon\gamma)/u}{m^2}.$$

We consider the following cases:

(i) If there exists $\alpha \in (0,1)$ for which $(1-\alpha)\varepsilon \geq m(\frac{\alpha\varepsilon}{u})$, then with $\beta = \alpha\varepsilon/u$ we have $m(\beta) \leq \varepsilon - \beta u$, which we already showed implies that the probability of interest is zero.

(ii) In view of (i), we assume without loss of generality that $(1-\alpha)\varepsilon < m(\frac{\alpha\varepsilon}{u})$, for all $\alpha \in (0,1)$. For any $\alpha \in (0,1)$, in the definition of $I(\varepsilon,\alpha,u)$ we may pick

$$\lambda = \frac{1}{1+\tau^2}\log\left(\frac{\tau^2 + x}{\tau^2(1-x)}\right),$$

which, after some algebra, yields

$$I(\varepsilon,\alpha,u) \geq H\left(\frac{x+\tau^2}{1+\tau^2}\bigg\|\frac{\tau^2}{1+\tau^2}\right),$$

where $H(y\|z) := y\log\frac{y}{z} + (1-y)\log\frac{1-y}{1-z}$ denotes the relative entropy between the Bernoulli($y$) and the Bernoulli($z$) distributions. This relative entropy is, in turn, by a standard argument (e.g., using Pinsker's inequality; cf. [2], Theorem 4.1), bounded below by $\frac{2x^2}{(1+\tau^2)^2}$. Therefore,

$$\frac{1}{n}\log\Pr\{S_n > n\varepsilon, |T_n| < nu\}$$

$$\leq -\sup_{\alpha\in(0,1)}\left[\frac{2x^2}{(1+\tau^2)^2}\right]$$

$$= -\sup_{\alpha\in(0,1)}\left[\frac{2(1-\alpha)^2\varepsilon^2}{m^2[1+(1+((\alpha\varepsilon)/u)^2 - (2\alpha\gamma\varepsilon)/u)/m^2]^2}\right],$$

proving part (ii).

(iii) Start by taking $u = K\varepsilon$. Noting that $|\gamma| \leq s \leq 1$,

$$1+\left(\frac{\alpha\varepsilon}{u}\right)^2 - \frac{2\gamma\alpha\varepsilon}{u} \leq 1+\left(\frac{\alpha\varepsilon}{u}\right)^2 + \frac{2\alpha\varepsilon}{u} = \left(1+\frac{\alpha}{K}\right)^2.$$

This and part (ii) show that the exponent of interest is bounded below by

$$2\sup_{\alpha\in(0,1)}\left[\frac{m\cdot(1-\alpha)}{m^2+(1+\alpha/K)^2}\right]^2\varepsilon^2,$$

where $m = m(\frac{\alpha\varepsilon}{u}) = m(\frac{\alpha}{K})$. Picking $\alpha = 1/2$, this is further bounded below by

$$\frac{1}{2}\left[\frac{M}{M^2 + (1 + 1/(2K))^2}\right]^2 \varepsilon^2,$$

where $M = m(\frac{1}{2K})$, giving the required result in the case $u = K\varepsilon$. Since the probability in (3.3) can be no bigger for smaller values of $u$, the same bound holds for all $0 < u \leq K\varepsilon$.  $\square$

We are now in a position to illustrate how the results of Proposition 1.1 stated in the introduction can be derived from Theorem 3.1.

PROOF OF PROPOSITION 1.1. Part (i) is already stated in (1.1), and part (ii) is immediate from Corollary 3.1. For parts (iii) and (iv) we will use the bound in Theorem 3.1(ii). To that end, we begin by defining two functions $F, U$ appropriately.

Recall that, for (iii), we only have the following information: $X$ is supported on $[1, \infty)$, $E(X) = 5/3$ and $\text{Var}(X) = 20/9$. Then we can define

$$U(x) := \frac{3x}{2\sqrt{5}} - \frac{\sqrt{5}}{2}, \qquad x \geq 1,$$

so that $E(U(X)) = 0$ and $\text{Var}(U(X)) = 1$. Noting that $\mu \geq 1$ and that $E[(X^{3/4})^2] \leq E(X^2) = \text{Var}(X) + E(X)^2 = 5$ implies that $\text{Var}(X^{3/4}) \leq 5 - 1 = 4$. Therefore, letting

$$F(x) := (x^{3/4} - \mu)/2, \qquad x \geq 1,$$

we have $E(F(X)) = 0$ and $\text{Var}(F(X)) \leq 1$. Using again the fact that $\mu \geq 1$, we obtain an upper bound on $m(\beta)$ as

$$m(\beta) \leq \sup_{x \geq 1}\left[\frac{x^{3/4}}{2} - \frac{1}{2} - \frac{3\beta x}{2\sqrt{5}} + \frac{\sqrt{5}\beta - 1}{2}\right].$$

This is a particularly easy maximization for $\beta \geq \frac{\sqrt{5}}{4}$, in which case the maximum is achieved at $x = 1$, giving

(3.6) $$m(\beta) \leq \tilde{m}(\beta) := \frac{\beta}{\sqrt{5}} \qquad \text{for } \beta \geq \frac{\sqrt{5}}{4}.$$

We can now apply (3.2). Let $S_n$ and $T_n$ be as in Theorem 3.1, and let $\hat{S}_n$ and $\hat{T}_n$ be as in the proposition. For arbitrary $\varepsilon > 0$ and $u = \frac{\varepsilon}{20}$, (3.2) gives

$$-\frac{1}{n}\log \Pr\left\{\hat{S}_n - \mu > \varepsilon, \left|\hat{T}_n - \frac{5}{3}\right| < u\right\}$$

$$(3.7) \quad = -\frac{1}{n}\log \Pr\{S_n > n\varepsilon/2, |T_n| < 3nu/2\sqrt{5}\}$$

$$\geq \frac{1}{2} \sup_{\alpha \in (0,1)} \left[\frac{(1-\alpha)\tilde{m}((20\sqrt{5}\alpha)/3)}{\tilde{m}((20\sqrt{5}\alpha)/3)^2 + 1 + ((20\sqrt{5}\alpha)/3)^2 - (40\sqrt{5}\alpha\gamma)/3}\right]^2 \varepsilon^2,$$

where $\gamma$ is the (yet unknown) covariance between $F(X)$ and $U(X)$. Restricting to $\alpha \geq 3/80$, using (3.6) and noting that $|\gamma| \leq 1$, the above exponent is further bounded below by,

$$\frac{1}{2} \sup_{3/80 \leq \alpha < 1} \left[\frac{20\alpha(1-\alpha)/3}{((20\alpha)/3)^2 + (1+(20\sqrt{5}\alpha)/3)^2}\right]^2 \varepsilon^2 \geq 0.005\varepsilon^2,$$

where the last inequality follows by taking $\alpha = 0.0552083$ in the above minimization (this $\alpha$ was selected by plotting the graph of the expression to be maximized and picking $\alpha$ to give a value near the maximum). This proves (iii) for $u = \varepsilon/20$, but, since the probability of interest is nondecreasing in $u$, the same bound holds for any $0 < u \leq \varepsilon/20$.

For part (iv), assuming that we also know that $\text{Cov}(X^{3/4}, X) = 20/21$, we can calculate

$$\gamma := \text{Cov}(F(X), U(X)) = \frac{3}{4\sqrt{5}} \text{Cov}(X^{3/4}, X) = \frac{\sqrt{5}}{7}.$$

From the bound in (3.7), restricting as before to $\alpha \geq 3/80$, using (3.6) and substituting the value of $\gamma$, gives

$$-\log \Pr\left\{\hat{S}_n - \mu > \varepsilon, \left|\hat{T}_n - \frac{5}{3}\right| < u\right\}$$

$$\geq \frac{n}{2} \sup_{3/80 \leq \alpha < 1} \left[\frac{20\alpha(1-\alpha)/3}{(2400\alpha^2)/9 - (200\alpha)/21 + 1}\right]^2 \varepsilon^2$$

$$\geq 0.0367n\varepsilon^2,$$

where the last inequality follows from choosing $\alpha = 0.0568$. This proves (iv) for $u = \varepsilon/20$, and, as before, the same bound remains valid for any $0 < u \leq \varepsilon/20$. □

**4. Bounds for light tails.** As before, let $S_n, T_n$ denote the partial sums of $\{F(X_i)\}, \{U(X_i)\}$, respectively, with respect to the i.i.d. random variables $X, X_1, X_2, \ldots$, with common law $P$. We assume that $E(F(X)) = E(U(X)) = 0$, and throughout this section we also assume that $F$ and $U$ have finite exponential moments, that is,

$$\Lambda(\theta) := \log E[e^{\theta F(X)}] < \infty,$$

and $E[e^{\theta U(X)}] < \infty$, for all $\theta \in \mathbb{R}$.

From Corollary 3.1 and the subsequent discussion, we know that the screened estimator always admits exponential error bounds,

$$(4.1) \quad \log \Pr\{S_n > n\varepsilon, |T_n| < nu\} \leq -n \max\{\Lambda_+^*(\varepsilon, u), \Gamma_+^*(\varepsilon, u)\}, \qquad n \geq 1,$$

for all $\varepsilon, u > 0$, where the exponents $\Lambda_+^*$ and $\Gamma_+^*$, given in Corollary 3.1 and Remark 1 after Corollary 3.1, respectively, are strictly positive. But in this setting, the standard estimates $\{\frac{1}{n}S_n\}$ also admit exponential error bounds; Cramér's theorem states that

$$(4.2) \qquad \log \Pr\{S_n > n\varepsilon\} \leq -n\Lambda^*(\varepsilon), \qquad n \geq 1,$$

where

$$\Lambda^*(\varepsilon) := \sup_{\theta \geq 0}\{\theta\varepsilon - \Lambda(\theta)\} > 0,$$

for any $\varepsilon > 0$; see [4]. Recall that the exponents in both (4.1) and (4.2) are asymptotically tight.

In this section we develop conditions under which the screened estimator offers a nontrivial improvement. That is, even when the error of the standard estimator decays exponentially, the error of the screened estimator has a better rate in the exponent. To that end, we look at difference

$$\Delta(\varepsilon, u) := \max\{\Lambda_+^*(\varepsilon, u), \Gamma_+^*(\varepsilon, u)\} - \Lambda^*(\varepsilon).$$

Clearly $\Delta(\varepsilon, u)$ is always nonnegative. Theorem 4.1 says that, as long as the covariance between $F(X)$ and $U(X)$ is nonzero, $\Delta(\varepsilon, u)$ is strictly positive for all $\varepsilon, u$ small enough. This is strengthened in Theorem 4.2, where it is shown that this improvement is a "first-order effect," in that, for small $\varepsilon, u$, $\Delta(\varepsilon, u)$ and $\max\{\Lambda_+^*(\varepsilon, u), \Gamma_+^*(\varepsilon, u)\}$ are each of order $\varepsilon^2$.

This leads to a different interpretation of the advantage offered by the screened estimator. Suppose that, for small $\varepsilon, u$, $\Lambda^*(\varepsilon) \approx c\varepsilon^2$, and that $\max\{\Lambda_+^*(\varepsilon, u), \Gamma_+^*(\varepsilon, u)\} \approx (c + c')\varepsilon^2$, for some $c, c' > 0$. Then for large $n$, the error of the standard estimator is

$$\Pr\{S_n > n\varepsilon\} \approx e^{-nc\varepsilon^2},$$

whereas for the screened estimator,

$$\Pr\{S_n > n\varepsilon, |T_n| < u\} \approx e^{-n(c+c')\varepsilon^2}.$$

In both cases, we have approximately Gaussian tails. Therefore, roughly speaking, we may interpret the result of Theorem 4.2 as saying that, as long as the covariance between $F(X)$ and $U(X)$ is nonzero, *the screened estimates are asymptotically Gaussian with a strictly smaller variance than the standard estimates.*

THEOREM 4.1. *Suppose that $E[F(X)] = E[U(X)] = 0$ and that $\gamma := \mathrm{Cov}(F(X), U(X))$ is nonzero. There exists $\varepsilon_0 > 0$ such that, for each $0 < \varepsilon < \varepsilon_0$, there exists $u_0 = u_0(\varepsilon) > 0$ such that $\Delta(\varepsilon, u) > 0$ for all $u \in (0, u_0)$.*

Note that the assumption on the covariance being nonzero cannot be relaxed. For example, let $X_i = Y_i Z_i$, $i \geq 1$, where $\{Y_i\}$ are i.i.d. nonnegative random variables, and $\{Z_i\}$ are i.i.d., independent of the $\{Y_i\}$, with each $Z_i = \pm 1$ with probability $1/2$. With $F(x) \equiv |x| - E|X_1|$ and $U(X) \equiv \mathrm{sign}(x)$, we have $F(X_i) = Y_i - E(Y_i)$ and $U(X_i) = Z_i$, so that $S_n$ and $T_n$ are independent for all $n \geq 1$. Therefore,

$$\Pr\{S_n > n\varepsilon, |T_n| < nu\} = \Pr\{S_n > n\varepsilon\} \Pr\{|T_n| < nu\},$$

and since $\lim_n \Pr\{|T_n| < nu\} = 1$, the exponents of the other two probabilities must be identical.

Whenever $\gamma$ is nonzero, the variances $\sigma^2(F)$, $\sigma^2(U)$ of $F(X)$ and $U(X)$, respectively, are both nonzero. If $\tilde{\Delta}(\varepsilon, u)$ denotes the corresponding difference of exponents for the normalized functions $F/\sigma(F)$ and $U/\sigma(U)$, then from the definitions,

$$\Delta(\varepsilon, u) = \tilde{\Delta}\left(\frac{\varepsilon}{\sigma(F)}, \frac{\varepsilon}{\sigma(U)}\right).$$

Therefore, in order to determine the nature of this difference for small $\varepsilon$ we can assume, without loss of generality, that $\mathrm{Var}(F(X)) = \mathrm{Var}(U(X)) = 1$.

THEOREM 4.2. *Suppose that $E[F(X)] = E[U(X)] = 0$, $\mathrm{Var}(F(X)) = \mathrm{Var}(U(X)) = 1$, and that $\gamma := \mathrm{Cov}(F(X), U(X))$ is nonzero. Then there exists $\alpha > 0$ such that*

$$\liminf_{\varepsilon \to 0} \frac{1}{\varepsilon^2} \Delta(\varepsilon, \alpha\varepsilon) > 0.$$

*In fact, there exists $\varepsilon_0 > 0$ such that*

$$\Delta\left(\varepsilon, \frac{|\gamma|}{4}\varepsilon\right) \geq \frac{\gamma^2}{8}\varepsilon^2,$$

*for all $\varepsilon \in (0, \varepsilon_0)$.*

Before giving the proofs of the theorems, we collect some technical facts in the following lemma.

LEMMA 4.1. *Suppose that $E[F(X)] = E[U(X)] = 0$ and that $\gamma := \mathrm{Cov}(F(X), U(X))$ is nonzero. Then:*

*(i) $\Lambda$ is smooth on $\mathbb{R}$, $\Lambda(0) = 0$, $\Lambda'(0) = 0$, $\lim_{\theta \to \infty} \Lambda'(\theta) = \bar{F} := \mathrm{ess\,sup}\, F(X)$, $\Lambda''(0) = \mathrm{Var}(F(X)) > 0$ and $\Lambda''(\theta) > 0$ for all $\theta \in \mathbb{R}$.*

(ii) *For each $0 < \varepsilon < \bar{F}$ there exists a unique $\theta^* = \theta^*(\varepsilon) > 0$ such that $\Lambda'(\theta^*) = \varepsilon$ and $\Lambda^*(\varepsilon) = \theta^*\varepsilon - \Lambda(\theta^*)$, where $\theta^* = \theta^*(\varepsilon)$ is strictly increasing in $\varepsilon \in (0, \bar{F})$.*

(iii) *Suppose $\mathrm{Var}(F(X)) = 1$. Let $\delta \geq 0$ be arbitrary (but fixed). Then for any $\eta > 0$ there exists $\bar{\varepsilon} > 0$ such that*

$$\Lambda(\delta\varepsilon) \geq \tfrac{1}{2}(1-\eta)\delta^2\varepsilon^2 \qquad \text{for all } \varepsilon < \bar{\varepsilon}.$$

(iv) *Suppose $\mathrm{Var}(F(X)) = \mathrm{Var}(U(X)) = 1$. For arbitrary (but fixed) $\beta \geq 0$, and for all $t, \varepsilon \geq 0$, define $f_t(\varepsilon) := \Lambda_+(t\varepsilon, \beta\varepsilon)$. Then for any $\eta > 0$ there exist $\tau, \bar{\varepsilon} > 0$ such that*

$$f_t(\varepsilon) \leq \tfrac{1}{2}(1 + \beta^2 - 2\beta\gamma + \eta)\varepsilon^2 \qquad \text{for all } \varepsilon < \bar{\varepsilon},\ |t-1| < \tau.$$

PROOF. The statements in (i) and (ii) are well known; see, for example, [4]. In particular, it is a standard exercise to apply the dominated convergence theorem in order to justify all the required differentiations, as well as all the continuity statements and differentiations in the rest of this proof and in the proofs of Theorems 4.1 and 4.2. For (iii), given $\eta > 0$, since $\Lambda''(\theta)$ is continuous and $\Lambda''(0) = 1$, we can choose $\varepsilon' > 0$ such that $\Lambda''(\varepsilon) \geq 1 - \eta$ for $\varepsilon < \varepsilon'$. The result follows upon expanding $\Lambda$ in a Taylor series around zero and recalling that $\Lambda(0) = \Lambda'(0) = 0$, with $\bar{\varepsilon} = \varepsilon'/\delta$.

Part (iv) is similar. Let $\eta > 0$ be given. We have $f_t(0) = \Lambda_+(0,0) = 0$, $f'_t(0) = E[tF(X) - \beta U(X)] = 0$ and $f''_t(\varepsilon)$ is jointly continuous in $t, \varepsilon \geq 0$ with

$$f''_t(0) = \mathrm{Var}(tF(X) - \beta U(X)) = t^2 + \beta^2 - 2t\beta\gamma,$$

where the prime ($'$) now denotes differentiation with respect to $\varepsilon$. Continuity at the point $(t, \varepsilon) = (1, 0)$ implies that we can find $\tau, \bar{\varepsilon} > 0$ such that

$$f''_t(\varepsilon) \leq f''_1(0) + \eta = 1 + \beta^2 - 2\beta\gamma + \eta \qquad \text{for all } \varepsilon < \bar{\varepsilon}, |t-1| < \tau.$$

For any $t$ in that range, expanding $f_t(\varepsilon)$ in a three-term Taylor series around $\varepsilon = 0$ gives the required result. $\square$

PROOF OF THEOREM 4.1. From the definitions, it follows that

$$(4.3) \quad \Delta(\varepsilon, u) \geq \Lambda^*_+(\varepsilon, u) - \Lambda^*(\varepsilon) \geq \sup_{\theta \geq 0}[-\theta u - \Lambda_+(\theta^*, \theta) + \Lambda(\theta^*)].$$

The expression inside the last supremum is zero for $\theta = 0$, and our goal is to show that it is strictly positive for small $\theta$. To that end, define the function

$$g(\theta) := E[F(X)U(X)e^{\theta F(X)}], \qquad \theta \geq 0,$$

and note that it is continuous in $\theta$, and $g(0) = \gamma$. Choose $\theta_0 > 0$ so that $g(\theta)/\gamma \geq 1/2$ for all $0 \leq \theta \leq \theta_0$. Let $\varepsilon_0 = \Lambda'(\theta_0) > 0$, and choose and fix an arbitrary $0 < \varepsilon < \varepsilon_0$, so that $\theta^* = \theta^*(\varepsilon) \in (0, \theta_0)$.

First consider the case $\gamma > 0$. Define
$$h(\theta) := \theta^*\varepsilon - \theta u - \Lambda_+(\theta^*, \theta).$$
Then $h(0) = \Lambda^*(\varepsilon)$, and as in (4.3),

$$\begin{aligned}
\Delta(\varepsilon, u) &\geq \Lambda_+^*(\varepsilon, u) - \Lambda^*(\varepsilon) \\
&\geq \left[\sup_{\theta \geq 0} h(\theta)\right] - \Lambda^*(\varepsilon) \\
&\geq h(0) - \Lambda^*(\varepsilon) = 0.
\end{aligned} \quad (4.4)$$

In order to establish that $\Delta(\varepsilon, u) > 0$ it suffices to show that $h'(0) > 0$. Computing the derivative of $h$ yields
$$h'(0) = e^{-\Lambda(\theta^*)} E[U(X)e^{\theta^* F(X)}] - u,$$
and expanding the exponential inside the expectation in a two-term Taylor expansion,
$$h'(0) = \theta^* e^{-\Lambda(\theta^*)} E[F(X)U(X)e^{\tilde\theta F(X)}] - u,$$
where $\tilde\theta = \tilde\theta(X) \in (0, \theta^*)$. Therefore,
$$h'(0) \geq \theta^* e^{-\Lambda(\theta^*)} \inf_{\theta \in (0,\theta^*)} g(\theta) - u \geq \theta^* e^{-\Lambda(\theta^*)} \gamma/2 - u,$$
which is strictly positive, as long as
$$u < u_0 = u_0(\varepsilon) := \theta^*(\varepsilon) e^{-\Lambda(\theta^*(\varepsilon))} |\gamma|/2.$$

The case $\gamma < 0$ is similar, with $\Gamma_+^*$ in place of $\Lambda_+^*$: Replace $h$ by $h(\theta) = \theta^*\varepsilon - \theta u - \log E[\exp\{\theta^* F(X) + \theta U(X)\}]$, so that $h(0) = \Lambda^*(\varepsilon)$ and
$$\Delta(\varepsilon, u) \geq \Gamma_+^*(\varepsilon, u) - \Lambda^*(\varepsilon) \geq \left[\sup_{\theta \geq 0} h(\theta)\right] - \Lambda^*(\varepsilon)$$
$$\geq h(0) - \Lambda^*(\varepsilon) = 0.$$
Again it suffices to show $h'(0) > 0$, where
$$\begin{aligned}
h'(0) &= -e^{-\Lambda(\theta^*)} E[U(X)e^{\theta^* F(X)}] - u \\
&= -\theta^* e^{-\Lambda(\theta^*)} E[F(X)U(X)e^{\tilde\theta F(X)}] - u,
\end{aligned}$$
with $\tilde\theta = \tilde\theta(X) \in (0, \theta^*)$. Then,
$$h'(0) \geq -\theta^* e^{-\Lambda(\theta^*)} \sup_{\theta \in (0,\theta^*)} g(\theta) - u \geq -\theta^* e^{-\Lambda(\theta^*)} \gamma/2 - u,$$
which is strictly positive, as long as $u < u_0 = u_0(\varepsilon)$, with the same $u_0$ as before. $\square$

PROOF OF THEOREM 4.2.  Assume first that $\gamma > 0$. Following the derivation of (4.4) in the proof of Theorem 4.1, we have that for any $0 < \varepsilon < \bar{F}$ and any $u, \phi > 0$,

$$(4.5) \qquad \Delta(\varepsilon, u) \geq -\phi u - \Lambda_+(\theta^*, \phi) + \Lambda(\theta^*).$$

At this point, most of the required work has been done. What remains is to write the above expression as a second-order Taylor expansion around $\varepsilon = 0$, so that, with $u = \gamma\varepsilon/4$ and $\phi = \gamma\varepsilon$, the right-hand side of (4.5) is approximately bounded below by

$$-\frac{\gamma^2\varepsilon^2}{4} - \frac{1}{2}\varepsilon^2\left[\frac{\partial^2}{\partial\varepsilon^2}\Lambda_+(\theta^*(\varepsilon), \gamma\varepsilon)\right]_{\varepsilon=0} + \frac{1}{2}\varepsilon^2\left[\frac{\partial^2}{\partial\varepsilon^2}\Lambda(\theta^*(\varepsilon))\right]_{\varepsilon=0} \geq \frac{\gamma^2\varepsilon^2}{8}.$$

We proceed to make this approximation rigorous. Let $\eta := \gamma^2/10 > 0$ in parts (iii) and (iv) of Lemma 4.1, and choose and fix a $\delta \in (0, \eta)$ smaller than the resulting $\tau$ in part (iv). Since $\Lambda''(\theta)$ is continuous and $\Lambda''(0) = 1$, we can choose $\theta_0 > 0$ small enough so that $|\Lambda''(\theta) - 1| \leq \delta\Lambda''(\theta)$ for all $0 < \theta < \theta_0$. Let $\varepsilon_0$ be the minimum of $\Lambda'(\theta_0)$ and the two quantities $\bar{\varepsilon}$ in parts (iii) and (iv) of the lemma. Then $\theta^*(\varepsilon) < \theta_0$ for all $0 < \varepsilon < \varepsilon_0$, and moreover, $\theta^*(\varepsilon) = \frac{\varepsilon}{\Lambda''(\theta)}$ for some $\theta < \theta^* < \theta_0$, so that

$$(4.6) \qquad \left|\frac{\theta^*(\varepsilon)}{\varepsilon} - 1\right| \leq \delta < \tau \qquad \text{for all } 0 < \varepsilon < \varepsilon_0.$$

Now for any $\varepsilon < \varepsilon_0$, let $u = \gamma\varepsilon/4$ and $\phi = \gamma\varepsilon$ in (4.5); using (4.6) and noting that $\Lambda(\theta^*)$ is nondecreasing in $\theta^*$,

$$\begin{aligned}
\Delta(\varepsilon, \gamma\varepsilon/4) &\geq -\frac{\gamma^2\varepsilon^2}{4} - \Lambda_+(\theta^*(\varepsilon), \gamma\varepsilon) + \Lambda((1-\delta)\varepsilon) \\
&\geq -\frac{\gamma^2\varepsilon^2}{4} - \frac{1}{2}(1 - \gamma^2 + \eta)\varepsilon^2 + \frac{1}{2}(1-\eta)(1-\delta)^2\varepsilon^2 \\
&= \frac{\varepsilon^2}{4}[\gamma^2 + 2(1-\eta)\delta^2 - 4(1-\eta)\delta - \eta] \\
&\geq \frac{\varepsilon^2}{4}[\gamma^2 - 5\eta] \geq \frac{\gamma^2\varepsilon^2}{8},
\end{aligned}$$

where the second inequality follows from parts (iii) and (iv) of Lemma 4.1 with $(1-\delta)$ in place of $\delta$, $\beta = \gamma$, and $t = \theta^*(\varepsilon)/\varepsilon$.

Finally, the same result holds in the case $\gamma < 0$, either by considering $-U$ in place of $U$, or by replacing $\Lambda_+^*$ by $\Gamma_+^*$ in the above argument, as in the proof of Theorem 4.1. □

**5. Proofs of Theorems 2.1, 2.2 and 2.3.** We begin with a simple, general upper bound in the spirit of the results in [3].

LEMMA 5.1. *Let $F_1, F_2, \ldots, F_m$ be an arbitrary (finite) collection of measurable functions from $\mathbb{R}$ to $\mathbb{R}$. For any constants $c_1, c_2, \ldots, c_m$ we have*

$$\log \Pr\left\{\sum_{i=1}^{n} F_j(X_i) > nc_j \text{ for all } j = 1, 2, \ldots, m\right\} \leq -n \inf_{Q \in E_m} H(Q\|P),$$

*where $E_m$ is the set of all probability measures $Q$ on $\mathbb{R}$ such that $\int F_j \, dQ > c_j$ for all $j = 1, 2, \ldots, m$.*

PROOF. Let $A$ denote the event of interest in the lemma, and assume without loss of generality that it has nonzero probability. Write $P_A$ for the probability measure on $\mathbb{R}^n$ obtained by conditioning the product measure $P^n$ on $A$, and note that, by definition,

$$-\log \Pr(A) = -\log P^n(A) = H(P_A \| P^n).$$

Expressing $P_A$ as the product of the conditional measures $P_{A,i}(\cdot | x_1, \ldots, x_{i-1})$ for $i = 1, 2, \ldots, n$, we can expand the logarithm inside the relative entropy to obtain

$$-\log \Pr(A) = \sum_{i=1}^{n} E[H(P_{A,i}(\cdot | Y_1, \ldots, Y_{i-1}) \| P)],$$

where the random variables $Y_1, Y_2, \ldots, Y_n$ have joint distribution given by the measure $P_A$. Using the fact that relative entropy is convex in its first argument (see, e.g., [4], Chapter 6), Jensen's inequality gives

$$-\log \Pr(A) \geq \sum_{i=1}^{n} H(Q_i \| P),$$

where $Q_i$ denotes the $i$th marginal of $P_A$ on $\mathbb{R}$. Using convexity again,

$$-\log \Pr(A) \geq n \sum_{i=1}^{n} \frac{1}{n} H(Q_i \| P) \geq n H(Q \| P),$$

where $Q = \frac{1}{n} \sum_{i=1}^{n} Q_i$. To complete the proof it suffices to show that $Q \in E_m$. Indeed, for any $j = 1, 2, \ldots, m$,

$$\int F_j \, dQ = \frac{1}{n} \sum_{i=1}^{n} \int F_j \, dQ_i = E\left[\frac{1}{n} \sum_{i=1}^{n} F_j(Y_i)\right] > c_j,$$

where the last inequality holds since the joint distribution of the $\{Y_i\}$ is $P_A$, which is entirely supported on $A$ by definition. $\square$

Next we give the proof of Theorem 2.2. The first upper bound follows from Lemma 5.1, the second is derived using the classical Chernoff bound, and the positivity of the exponent comes from the domination assumption $m(\beta) < \infty$.

PROOF OF THEOREM 2.2. Throughout, we assume, without loss of generality, that $\mu = \nu = 0$. For part (i), taking $F_1 = F$, $F_2 = -U$, $c_1 = \varepsilon$ and $c_2 = -u$, Lemma 5.1 immediately yields the required bound. Part (ii) follows by the usual Chernoff argument: For any pair of $\theta_1, \theta_2 \geq 0$,

$$\begin{aligned}
\Pr\{S_n > n\varepsilon, T_n < nu\} &\leq \Pr\{S_n > n\varepsilon, T_n < nu\} \\
&= E[\mathbb{I}_{\{S_n > n\varepsilon\}}\mathbb{I}_{\{T_n < nu\}}] \\
&\leq E[\exp\{\theta_1(S_n - n\varepsilon)\}\exp\{-\theta_2(T_n - nu)\}] \\
&= \exp\{-n[\theta_1\varepsilon - \theta_2 u - \Lambda_+(\theta_1, \theta_2)]\}.
\end{aligned}$$

The stated result is obtained upon taking the supremum over all $\theta_1, \theta_2 \geq 0$ in the exponent.

Finally, for part (iii) choose and fix an arbitrary $\alpha \in (0,1)$. Taking $\theta_2 = \alpha\varepsilon\theta_1/u$ in the definition of $\Lambda_+^*(\varepsilon, u)$ yields

(5.1) $$\Lambda_+^*(\varepsilon, u) \geq \sup_{\theta \geq 0}[\theta(1-\alpha)\varepsilon - \Lambda_0(\theta)],$$

where $\Lambda_0(\theta) := \Lambda_+(\theta, \frac{\alpha\varepsilon\theta}{u}) < \infty$ for all $\theta \geq 0$ because $m(\beta) < \infty$ for all $\beta > 0$. Now for any $\theta \geq 0$, let $X_\theta$ be a random variable whose distribution has Radon–Nikodym derivative with respect to that of $X$ given by the density

$$g_\theta(x) = \frac{\exp\{\theta[F(x) - ((\alpha\varepsilon)/u)U(x)]\}}{E[\exp\{\theta[F(X) - ((\alpha\varepsilon)/u)U(X)]\}]}, \qquad x \in \mathbb{R},$$

so that $g_0 \equiv 1$ and $X_0 = X$. Obviously $\Lambda_0(0) = 0$, and simple calculus shows that $\Lambda_0'(\theta) = E[F(X_\theta) - \frac{\alpha\varepsilon}{u}U(X_\theta)]$ so that $\Lambda_0'(0) = 0$; the dominated convergence theorem justifies the differentiation under the integral, and also shows that $\Lambda_0'(\theta)$ is continuous in $\theta$ for all $\theta \geq 0$, since $F(X)$ and $U(X)$ have finite first moments and $m(\beta) < \infty$ for all $\beta > 0$.

Pick $\theta_0 > 0$ small enough so that

$$\sup\{\Lambda_0'(\theta) : \theta \in [0, \theta_0]\} \leq \Lambda_0'(0) + \frac{(1-\alpha)\varepsilon}{2}.$$

Restricting the range of the supremum in (5.1) to $[0, \theta_0]$ yields

$$\Lambda_+^*(\varepsilon, u) \geq \sup_{0 \leq \theta \leq \theta_0} \frac{\theta(1-\alpha)\varepsilon}{2} = \frac{\theta_0(1-\alpha)\varepsilon}{2},$$

which is strictly positive. □

The main technical step in the following proof is the (asymptotic) large deviations lower bound; it is established by a change-of-measure argument combined with regularization of the random variables of interest, as in Cramér's theorem. The main difference from the classical case is that, here, the domination assumption $m(\beta) < \infty$ replaces the usual condition on the existence of exponential moments in a neighborhood of the origin.

PROOF OF THEOREM 2.3. As above, we assume without loss of generality that $\mu = \nu = 0$. Write $\boldsymbol{\theta}$ for an arbitrary pair of nonnegative $(\theta_1, \theta_2)$, and write $\mathbf{G}: \mathbb{R} \to \mathbb{R}^2$ for the function $\mathbf{G}(x) = (F(x), -U(x))$, $x \in \mathbb{R}$, so that $\Lambda_+(\boldsymbol{\theta}) = \log E[\exp\{\langle \boldsymbol{\theta}, \mathbf{G}(X) \rangle\}]$ and

$$\Lambda_+^*(\varepsilon, u) = \sup_{\boldsymbol{\theta}} [\langle \boldsymbol{\theta}, (\varepsilon, -u) \rangle - \Lambda_+(\boldsymbol{\theta})],$$

where $\langle \cdot, \cdot \rangle$ denotes the usual Euclidean inner product. Note that, since $m(\beta) < \infty$ for all $\beta > 0$, we have $\Lambda_+(\boldsymbol{\theta}) < \infty$ as long as $\theta_2 > 0$, and $\Lambda_+(\mathbf{0}) = 0$. Moreover, since $E(\mathbf{G}(X)) = \mathbf{0}$, the dominated convergence theorem implies that $\Lambda_+(\boldsymbol{\theta})$ is differentiable, with

(5.2) $$\nabla \Lambda_+(\boldsymbol{\theta}) = E[\mathbf{G}(X) \exp\{\langle \boldsymbol{\theta}, \mathbf{G}(X) \rangle - \Lambda_+(\boldsymbol{\theta})\}],$$

for all $\boldsymbol{\theta}$ with $\theta_2 > 0$.

In view of Theorem 2.2(ii), in order to establish the limiting relation (2.2), it suffices to prove the asymptotic lower bound,

(5.3) $$\liminf_{n \to \infty} \frac{1}{n} \log \Pr\{S_n > n\varepsilon, T_n < nu\} \geq -\Lambda_+^*(\varepsilon, u).$$

To that end, consider three cases. First, if $\Lambda_+^*(\varepsilon, u) = \infty$, (5.3) is trivially true. Second, assume that $\Lambda_+^*(\varepsilon, u) < \infty$ and there exists $\boldsymbol{\theta}$ such that

(5.4) $$E[\mathbf{G}(X) \exp\{\langle \boldsymbol{\theta}, \mathbf{G}(X) \rangle - \Lambda_+(\boldsymbol{\theta})\}] = (\varepsilon, -u).$$

Fixing this $\boldsymbol{\theta}$, define a new sequence of i.i.d. random variables $X', X_1', X_2', \ldots$ with common distribution $P'$, where

$$\frac{dP'}{dP}(x) = \exp\{\langle \boldsymbol{\theta}, \mathbf{G}(x) \rangle - \Lambda_+(\boldsymbol{\theta})\}, \qquad x \in \mathbb{R}.$$

Write $S_n'$ and $T_n'$ for the corresponding partial sums, and choose and fix $\delta > 0$; then, $\frac{1}{n} \log \Pr\{S_n > n\varepsilon, T_n < nu\}$ is bounded below by

(5.5)
$$\frac{1}{n} \log \Pr\{n\varepsilon < S_n < n(\varepsilon + \delta), n(u - \delta) < T_n < nu\}$$
$$= \frac{1}{n} \log E\left[\prod_{i=1}^n \frac{dP}{dP'}(X_i') \mathbb{I}_{\{n\varepsilon < S_n' < n(\varepsilon + \delta)\}} \mathbb{I}_{\{n(u-\delta) < T_n' < nu\}}\right]$$
$$= \Lambda_+(\boldsymbol{\theta}) - \langle \boldsymbol{\theta}, (\varepsilon, -u) \rangle + \frac{1}{n} \log E[e^{-\theta_1(S_n' - n\varepsilon) + \theta_2(T_n' - nu)} \mathbb{I}_{B_n}]$$
$$\geq \Lambda_+(\boldsymbol{\theta}) - \langle \boldsymbol{\theta}, (\varepsilon, -u) \rangle - (\theta_1 + \theta_2)\delta + \frac{1}{n} \log \Pr(B_n),$$

where $B_n$ denotes the event $B_n := \{n\varepsilon < S'_n < n(\varepsilon + \delta)\} \cap \{n(u - \delta) < T'_n < nu\}$, and the last inequality follows from the observation that the exponential inside the expectation is bounded below by $\exp\{-\theta_1 n\delta - \theta_2 n\delta\}$ on $B_n$. Note that our assumption (5.4) implies that $E[\mathbf{G}(X')] = (\varepsilon, -u)$, and since $m(\beta) < \infty$ for all $\beta$, $F(X')$ and $U(X')$ have finite second moments. Therefore, from the central limit theorem we obtain

$$\liminf_{n\to\infty} \frac{1}{n} \log \Pr(B_n) = 0,$$

as long as $\delta > 0$ is fixed. Noting also that $\Lambda_+(\boldsymbol{\theta}) - \langle \boldsymbol{\theta}, (\varepsilon, -u) \rangle \geq -\Lambda_+^*(\varepsilon, u)$, taking $n \to \infty$ in (5.5) we obtain

$$(5.6) \quad \liminf_{n\to\infty} \frac{1}{n} \log \Pr\{S_n > n\varepsilon, T_n < nu\} \geq -\Lambda_+^*(\varepsilon, u) - (\theta_1 + \theta_2)\delta,$$

for each $\delta > 0$, and taking $\delta \downarrow 0$ in the above right-hand side yields (5.3).

The third and last case is when $\Lambda_+^*(\varepsilon, u) < \infty$ but there is no $\boldsymbol{\theta}$ such that (5.4) is satisfied. We will repeat the above argument, but instead of the sequence $\{\mathbf{G}(X_n)\}$ we will consider the new i.i.d. sequence $\{\mathbf{H}(X_n)\}$ which is obtained by adding to the $\{\mathbf{G}(X_n)\}$ i.i.d. Gaussians with small mean and variance. Specifically, choose and fix arbitrary $\delta > 0$ and $t > 0$, and let

$$\mathbf{H}(X_n) := \mathbf{G}(X_n) + t\mathbf{Z}_n + \left(\frac{\delta}{2}, \frac{\delta}{2}\right), \qquad n \geq 1,$$

where the $\{\mathbf{Z}_n\}$ are i.i.d. with each $\mathbf{Z}_n$ consisting of two independent standard Gaussian components, independent of the $\{X_n\}$. Let

$$\Lambda_t(\boldsymbol{\theta}) := \log E[\exp\{\langle \boldsymbol{\theta}, \mathbf{H}(X)\rangle\}],$$

and note that

$$(5.7) \quad \Lambda_t(\boldsymbol{\theta}) = \Lambda_+(\boldsymbol{\theta}) + t^2(\theta_1^2 + \theta_2^2)/2 + \delta(\theta_1 + \theta_2)/2 \geq \Lambda_+(\boldsymbol{\theta}) \geq 0,$$

where the last inequality follows by applying Jensen's inequality to the logarithm in the definition of $\Lambda_+(\boldsymbol{\theta})$ and recalling that $\mathbf{G}(X)$ has zero mean. Consequently,

$$(5.8) \qquad \Lambda_t^*(\varepsilon, u) := \sup_{\boldsymbol{\theta}}[\langle \boldsymbol{\theta}, (\varepsilon, -u)\rangle - \Lambda_t(\boldsymbol{\theta})] \leq \Lambda_+^*(\varepsilon, u) < \infty.$$

From (5.7) and (5.8) it follows that, for any given $\boldsymbol{\theta}$, the function

$$L(\boldsymbol{\theta}) := \langle \boldsymbol{\theta}, (\varepsilon, -u)\rangle - \Lambda_t(\boldsymbol{\theta}) \leq \Lambda_+^*(\varepsilon, u) - t^2(\theta_1^2 + \theta_2^2)/2 - \delta(\theta_1 + \theta_2)/2$$

has $\sup_{\boldsymbol{\theta}: \theta_1 + \theta_2 > R} L(\boldsymbol{\theta}) \to -\infty$ as $R \to \infty$. Moreover, in view of (5.2), $L(\boldsymbol{\theta})$ is differentiable, and therefore the supremum in the definition of $\Lambda_t^*(\varepsilon, u)$ is achieved for some finite $\boldsymbol{\theta}$ which satisfies the analog of (5.4), that is, with $\mathbf{H}$ and $\Lambda_t(\boldsymbol{\theta})$ in place of $\mathbf{G}$ and $\Lambda_+(\boldsymbol{\theta})$, respectively. So we can conclude from the previous argument that the lower bound (5.3) holds with $\mathbf{H}$ in place of **G**. In fact, for the specific value of $\delta > 0$ we chose in the definition of **H**, the same argument used to establish (5.5) and then (5.6) yields the following asymptotic lower bound:

$$\liminf_{n\to\infty} \frac{1}{n} \log \Pr\left\{ n\varepsilon < S_n + t\sqrt{n}W + \frac{n\delta}{2} < n(\varepsilon + \delta), \right.$$

$$\left. n(u - \delta) < T_n + t\sqrt{n}V - \frac{n\delta}{2} < nu \right\}$$

(5.9)
$$\geq -\Lambda_t^*(\varepsilon, u) - (\theta_1 + \theta_2)\delta$$
$$\geq -\Lambda_+^*(\varepsilon, u) - (\theta_1 + \theta_2)\delta$$
$$> -\infty,$$

where $W, V$ are independent standard Gaussian random variables independent of the $\{X_n\}$. On the other hand, a simple union bound gives

$$\Pr\left\{ n\varepsilon < S_n + t\sqrt{n}W + \frac{n\delta}{2} < n(\varepsilon + \delta), \right.$$

$$\left. n(u - \delta) < T_n + t\sqrt{n}V - \frac{n\delta}{2} < nu \right\}$$

(5.10)
$$\leq \Pr\{n\varepsilon < S_n < n(\varepsilon + 2\delta), n(u - 2\delta) < T_n < nu\}$$
$$+ \Pr\left\{ |W| \geq \frac{\sqrt{n}\delta}{2t}, |V| \geq \frac{\sqrt{n}\delta}{2t} \right\},$$

where the last probability is easily bounded as

(5.11) $$\frac{1}{n} \log \Pr\left\{ |W| \geq \frac{\sqrt{n}\delta}{2t}, |V| \geq \frac{\sqrt{n}\delta}{2t} \right\} \leq -\frac{\delta^2}{4t^2}.$$

Combining the bounds (5.9), (5.10) and (5.11) yields

$$-\Lambda_+^*(\varepsilon, u) - (\theta_1 + \theta_2)\delta$$
$$\leq \max\left\{ -\frac{\delta^2}{4t^2}, \right.$$
$$\left. \liminf_{n\to\infty} \frac{1}{n} \log \Pr\{n\varepsilon < S_n < n(\varepsilon + 2\delta), n(u - 2\delta) < T_n < nu\} \right\}.$$

Letting $t \downarrow 0$ implies that

$$\liminf_{n\to\infty} \frac{1}{n} \log \Pr\{S_n > n\varepsilon, T_n < nu\} \geq -\Lambda_+^*(\varepsilon, u) - (\theta_1 + \theta_2)\delta,$$

and letting $\delta \downarrow 0$ establishes (5.3) and thus completes the proof of (2.2).

Finally, in order to show that the two rate functions are identical, it suffices to show that $\Lambda_+^*(\varepsilon, u)$ is no greater than the entropy $H(E\|Q)$, since the reverse inequality follows from the upper bound in Theorem 2.2(i) combined with the asymptotic relation (2.2) we just established. Indeed, for arbitrary $\theta_1, \theta_2 \geq 0$ and any $Q \in E$,

$$\theta_1 \varepsilon - \theta_2 u - \log E[\exp\{\theta_1 F(X) - \theta_2 U(X)\}]$$
$$= \theta_1 \varepsilon - \theta_2 u - \log \int dQ(x) \frac{dP}{dQ}(x) \exp\{\theta_1 F(x) - \theta_2 U(x)\}$$
$$\leq \theta_1 \varepsilon - \theta_2 u - \int dQ(x) \log\left[\frac{dP}{dQ}(x) \exp\{\theta_1 F(x) - \theta_2 U(x)\}\right]$$
$$= \theta_1\left[\varepsilon - \int F\, dQ\right] - \theta_2\left[u - \int U\, dQ\right] + H(Q\|P)$$
$$\leq H(Q\|P),$$

where the first inequality is simply Jensen's inequality and the second follows from the assumption that $Q \in E$. Taking the supremum of both sides over all $\theta_1, \theta_2 \geq 0$ and then the infimum over all $Q \in E$ establishes the inequality $\Lambda_+^*(\varepsilon, u) \leq H(E\|P)$ and completes the proof. □

It is now a simple matter to deduce Theorem 2.1 from Theorems 2.2 and 2.3.

PROOF OF THEOREM 2.1. Again we assume without loss of generality that $\mu = \nu = 0$. For part (i), since $E[e^{\theta F(X)}]$ is infinite for all $\theta > 0$, it is well known that

(5.12) $$\lim_{n\to\infty} \frac{1}{n} \log \Pr\{S_n > n\varepsilon\} = 0:$$

see, for example, [5], Example 9.8, page 78. To see that $H(\Sigma\|P) := \inf_{Q\in\Sigma} H(Q\|P) = 0$ note that, from Lemma 5.1, we have $\log \Pr\{S_n > n\varepsilon\} \leq -nH(\Sigma\|P)$. This combined with (5.12) implies that $H(\Sigma\|P) = 0$. The limit in part (ii) is an immediate consequence of Theorem 2.3, and the fact that the exponent is strictly nonzero follows from Theorem 2.2(iii) and the identification of the rate function as the entropy given in Theorem 2.3. □

**Acknowledgment.** We thank Peter Glynn and Jose Blanchet for several interesting conversations.

DEPARTMENT OF INFORMATICS
ATHENS UNIVERSITY OF ECONOMICS AND BUSINESS
PATISSION 76, ATHENS 10434
GREECE
E-MAIL: yiannis@aueb.gr
URL: http://pages.cs.aueb.gr/users/yiannisk/

DEPARTMENT OF ELECTRICAL
AND COMPUTER ENGINEERING
AND THE COORDINATED SCIENCES LABORATORY
UNIVERSITY OF ILLINOIS,
URBANA, ILLINOIS 61801
USA
E-MAIL: meyn@uiuc.edu
URL: http://black.csl.uiuc.edu/~meyn/